\magnification=\magstep1
\baselineskip 20pt
\input amstex 
\documentstyle{amsppt}

\def\sqep{\qed}

\topmatter
\title
The index of operators on foliated bundles
\endtitle

\author Victor Nistor \endauthor

\address
Department of Mathematics, Pennsylvania State University,
University Park, PA 16802 
\endaddress
\thanks { } Partially supported by NSF grant DMS-9205548,
a Sloan research fellowship
and NSF Young Investigator Award DMS-9457859
\endthanks

\abstract
We compute the equivariant cohomology Chern character
of the index of elliptic operators along the leaves of the 
foliation of a flat bundle. The proof is based on the study
of certain algebras of pseudodifferential operators and uses 
techniques for analizing noncommutative algebras similar to 
those developed in Algebraic Topology, but in the framework of
cyclic cohomology and noncommutative geometry. 
\endabstract
\endtopmatter

\document
\baselineskip = 20pt 
\TagsOnRight


\subheading{1. Introduction}

\bigpagebreak

Let $(V, \pi)$,\ $F \to V\ 
{\buildrel\pi\over\longrightarrow}\ B $ be a smooth
fiber bundle  with fiber $F$ of dimension 
$q$.  We assume that  $(V,\pi)$  is endowed with a flat connection 
corresponding to an integrable subbundle  ${\Cal F} \subset TV$, 
of dimension 
$n = dim(B)$, transverse at any point to the fibers of $\pi$.  
The pair $(V,{\Cal F})$  is a foliation. 

The purpose of this paper to study invariants of differential operators
along the leaves  of the above foliation.
The index of an elliptic operator along the leaves of the foliation
${\Cal F}$ is an element in the $K$-Theory group
$K^0({\Cal F})=K_0(\Psi^{-\infty}({\Cal F}))$ where 
$\Psi^{-\infty}({\Cal F})$ is the algebra of regularizing 
operators along the leaves. In the case of a foliated
bundle there exists a Connes-Karoubi character
$Ch:K_0(\Psi^{-\infty}({\Cal F}))\to H_\Gamma^{*-q}(F,{\Cal O})'$ to 
the dual of equivariant
cohomology with twisted coefficients, where $\Gamma$
is the fundamental group of the base $B$ acting on the fiber $F$
via holonomy, and ${\Cal O}$ is the orientation sheaf. 
Our main theorem computes
the Connes-Karoubi character of the index.
This amounts to a proof of the ``higher 
index theorem for foliations'' in this special case.
A very general higher index theorem for foliations can be found in \cite{C3}
and here we give a new proof of this theorem for flat bundles.
Some very interesting results related to the results in this paper 
are contained in \cite{CM2} where $Diff$-invariant structures are
treated in detail. See also \cite{C4}.
In the even more special case of a family of elliptic
operators our theorem recovers the computation of the Chern character
of the index bundle of a family of elliptic operators
\cite{AS}.

The problem that we consider was suggested by 
\cite{DHK}. 
The proof of our theorem is based on the Cuntz-Quillen exact 
sequence \cite{CQ} and the results in \cite{BN} and 
\cite{N3}.  The second  paper identifies the cyclic cohomology groups with 
geometric groups.  The third paper provides us with the axiomatic setting 
necessary to deal with index problems in the framework of cyclic cohomology.

Let $P(\partial)=\sum a_\alpha \partial^\alpha$ be the local expression of
an elliptic operator on the base $B$, acting between the sections of
two vector bundles. We lift each vector field $\partial_i$ on $B$
to a vector field $\nabla_i$ on the total space $V$ of our flat bundle.
This will allow us to construct the lift $P(\nabla)$ which will be an
example of an elliptic differential operator along the leaves 
of ${\Cal F}$. In this rather degenerate case the invariants for 
$P(\nabla)$ reduce to invariants of $P(\partial)$. However not all
operators that we consider arise in this way, actually very few do.
The nonmultiplicativity of the signature \cite{A}
is related to the fenomena that we investigate.

\bigpagebreak

\subheading{2. Statement of the problem}

\bigpagebreak

Consider a smooth foliation ${\Cal F}$  of smooth manifold $V$.  All the 
structures that will be used in this paper will be smooth, i.e.
$C^\infty$, so that we shall omit ``smooth'' in the following.  
We think of the 
foliation $(V, {\Cal F})$  as an integrable subbundle ${\Cal F} \subset 
TV$.  That is, ${\Cal F}$  identifies with the tangent bundle to the 
foliation.

By considering only differentiations along the fibers of ${\Cal F}$  one 
obtains longitudinal differential operators.  In analogy with manifolds,
one can proceed then to define longitudinal pseudodifferential operators, 
denoted  $\Psi^p({\Cal F})$. A good reference to 
these constructions is  \cite{MoS}. An alternative description of
these algebras for foliations coming from flat bundles is given in
the next section. 

The algebra  $\Psi^{-\infty}({\Cal F})$  of  regularizing operators along 
${\Cal F}$  is usually referred to as the ``algebra of the foliation''.  
If, by abuse of notation, we use the symbol ${\Cal F}$  to denote the graph 
of the foliation  $(V, {\Cal F})$  as well  \cite{W} then 
$\Psi^{-\infty}({\Cal F})$  identifies with the algebra $C_c^\infty({\Cal F})$ 
of compactly supported smooth kernels on the graph.  

We review here, in order to fix notation, the construction of the graph 
of the foliation  $(V, {\Cal F})$.  It consists of equivalence classes of 
triples  $(x,y,\gamma)$, where  $x,y \in V$  are on the same leaf and $\gamma$ 
is a path from $x$ to $y$ completely contained in that leaf.  The equivalence 
is given by ``holonomy''.  The graph is a smooth manifold, usually non Hausdorff.

As in the classical case, the principal symbol induces an isomorphism
$$
\sigma_0:  \Psi^0({\Cal F})/\Psi^{-1}({\Cal F}) \longrightarrow 
C_c^\infty(S^*{\Cal F})
$$
where  $S^*{\Cal F}$  is the unit sphere bundle in the dual ${\Cal F}^*$ of 
${\Cal F}$.  The notion of asymptotic expansion generalizes as well, and this 
shows that a matrix of order 0 pseudodifferential operators is invertible 
modulo regularizing operators if and only if its principal symbol
is invertible. 
From this we infer that  $\sigma_0$  induces an isomorphism 
$$
\sigma_{0*}: K_1^{\text{top}}(\Psi^0({\Cal F})/\Psi^{-\infty}({\Cal F})) 
\longrightarrow K_1^{\text{top}}(C_c^\infty(S^*{\Cal F}))
$$
where  $K_i^{\text{top}}$  is the quotient of  $K_i^{\text{alg}}$  
with respect 
to homotopy ($i = 0$ or $i = 1$).  

The most general form of the index problem for foliation is:

\bigpagebreak
\noindent
{\bf (FOL-ALG)}
\quad {\it Determine the algebraic $K$-theory boundary (index) map}
$$
Ind:  K_1^{\text{alg}}(\Psi^0({\Cal F})/\Psi^{-\infty}({\Cal F}))\ 
\longrightarrow\ K_0^{\text{alg}}(\Psi^{-\infty}({\Cal F})) .
$$

(A very closed related problem is obtained by considering topological 
$K$-theory.)

\bigpagebreak

The major difficulty is that we know very 
little about the $K_0$-groups involved. (The $K_1$-groups are 
relatively easy to determine.)

Denote by $HC_i^{per}(A)$, $i \in {\Bbb Z}/2{\Bbb Z}$, 
the periodic cyclic 
homology groups of an arbitrary complex algebra $A$, and by
$$
Ch: K_i^{alg}(A) \longrightarrow HC_i^{per}(A)
$$
the Connes-Karoubi character \cite{C2,K,LQ}.
One way to avoid the above difficulty is to ask 
for the Connes-Karoubi character of the index. 

\bigpagebreak  

\noindent
{\bf  (FOL-COH)}\quad {\it Compute the composition}
$$
Ch \circ Ind:  K_1^{\text{top}} (C_c^\infty(S^*{\Cal F})) \longrightarrow 
HC_0^{per}(\Psi^{-\infty}({\Cal F})) .
$$

\bigpagebreak

A formula for $Ch \circ Ind$  will be called a ``cohomological index theorem''.

The cohomological form of the problem is not just 
a simplification of the original 
problem, but it also brings a new perspective.  
This is because what we usually 
want is   $Ch(Ind [u])$, not  $Ind[u]$ itself. Also this
form of the problem 
makes the connection with the characteristic classes 
of foliations, as we shall see bellow.

The actual definition of the various cyclic homology groups will 
{\it not} be necessary 
for our purposes.  What will matter will be that they exist, and that they 
satisfy certain general properties.   This is very similar to the philosophy 
of Algebraic Topology, especially in the axiomatic 
approach due to Eilenberg and Steenrod.  

Let us begin by explaining some of the constructions in a particular but
suggestive case. Let $A$ be an algebra and let
$\tau: A \to {\Bbb C}$  be a trace (i.e., $\tau(xy) = 
\tau(yx))$.  The map that associates to any idempotent  $e = (e_{ij}) \in 
M_n(A)$  the number  $\tau(e) = \Sigma \tau(e_{ii}) 
\in {\Bbb C}$  factors to 
a morphism
$$
\tau_*:  K_0^{alg}(A) \longrightarrow {\Bbb C}.
$$

In general any trace $\tau$  defines a class 
$[\tau] \in HC^0_{per}(A)$  and there exists
a pairing $\langle\ ,\ \rangle$  between cycle homology and cyclic 
cohomology such that  
$$
\tau_*(e) = \langle Ch([e]), [\tau]\rangle
$$
We shall call the elements of $HC^i_{per}(A)$ higher traces.

Any holonomy invariant measure $\mu$  on a foliation 
$(V, {\Cal F})$  determines 
a trace  $\tau_\mu: \Psi^{-\infty}({\Cal F}) \to {\Bbb C}$.  The quantity 
$\tau_\mu(Ind[u])$  was determined by Connes in \cite{C1}.  In view of what we 
said above this amounts to a partial determination of  
$Ch(Ind[u])$.

\bigpagebreak

We now very briefly review the most important properties of periodic cyclic 
(co)homology to be used in the following.

\item{(1)}  The groups  $HC_i^{per}(A)$  and  $HC^i_{per}(A)$, 
$i \in {\Bbb Z}_2$ 
are covariant (resp. contravariant) functors on the category of 
complex locally 
convex algebras with continuous algebra morphisms.
If $f:A\to B$ is an algebra morphism then we denote by 
$f_*$ and, respectively, $f^*$ the induced morphisms.  

\item{(2)}  There is a pairing  $\langle -,-\rangle: HC_i^{per}(A) \otimes 
HC^i_{per}(A) \to {\Bbb C}$.

\item{(3)}  
$$
HC_i^{per}(\Bbb C) \simeq
\cases {\Bbb C} &\qquad  i = 0 \\
0 & \qquad i = 1
\endcases
$$
such that  $HC^0_{per}({\Bbb C})$  is generated by the identity map (trace).
$$
HC_i^{per}({\Bbb C}[{\Bbb Z}]) \simeq {\Bbb C}\qquad \text{for\ any}\quad 
i \in {\Bbb Z}/2{\Bbb Z},
$$
where for any group $\Gamma$  we denote by  ${\Bbb C}[\Gamma]$  its group 
algebra.

\item{(4)}  Consider a separated \'etale groupoid  ${\Cal G}$ \cite{BN}, that is, 
${\Cal G}$  is a small category together with manifold structures on 
${\Cal G}^{(0)}\ {\dsize\buildrel{\text{def}}\over =}\ Ob({\Cal G})$  and 
${\Cal G}^{(1)}\ {\dsize\buildrel{\text{def}}\over =}\ Mor({\Cal G})$  such 
that all morphisms are invertible, all structural maps are smooth and the 
domain and range are local diffeomorphisms.  Let  $B{\Cal G}$  be the 
geometric realization of the nerve of  ${\Cal G}$  (this is Grothendieck's 
classifying space of the topological category  ${\Cal G}$).  Also denote by 
${\Cal O}({\Cal G})$  the complex orientation sheaf of  $B{\Cal G}$  
(this is defined 
because  ${\Cal G}$  is \'etale) and denote by  $q$  the common 
dimension of  ${\Cal G}^{(0)}$ and ${\Cal G}^{(1)}$.  The main result of 
\cite{BN} establishes the existence of an injective map
$$
\Phi: {\bigoplus_{i+j \equiv q(2)}} H^j(B{\Cal G},{\Cal O}({\Cal G})) 
\longrightarrow HC^i_{per}(C_c^\infty({\Cal G}))
$$
where  $C_c^\infty({\Cal G})$  is endowed with the natural topology and the 
convolution product.  The morphism $\Phi$  is multiplicative and functorial 
with respect to {\it \'etale} morphisms \cite{N3}
that are one-to-one on units.
We will call $\Phi$ {\it the geometric map}.
The morphism $\Phi$ is multiplicative, $\Phi(\xi \times \zeta)=
\Phi(\xi)\otimes \Phi(\zeta)$ where $\times$ is the external
product in cohomology and $\otimes$ is the external product in 
periodic cyclic cohomology \cite{N2,N3}.

\item{(4')}  We are going to make more explicit the constructions 
of (4) in the 
case of interest for us.  Let $\Gamma$  be a discrete group acting on a manifold 
$X$.  Define  ${\Cal G}$  by  ${\Cal G}^{(0)} = X$,\ ${\Cal G}^{(1)} = 
X \times \Gamma$  with $domain(x,\gamma) = x$, $range(x,\gamma) = \gamma 
\cdot x$.  Then  $C_c^\infty({\Cal G})$  is the (algebraic) crossed product 
algebra and  $B{\Cal G} = (X \times E\Gamma)/\Gamma =$ ``the homotopy quotient 
$X//\Gamma$'', here  $\Gamma \to E\Gamma \to B\Gamma$  
is the universal principal 
$\Gamma$-bundle.  The map $\Phi$ of (4) becomes an injective map
$$
\Phi: {\bigoplus_{i+j \equiv  dim(X)(2)}} H^j(X//\Gamma, 
{\Cal O}(X)) \to 
HC^i_{per}(C^\infty(X)\rtimes \Gamma)
$$
(this is a more precise form of some results in \cite{N1}).  If $X$  is orientable 
and $\Gamma$ preserves the orientation, then the left hand side reduces to 
equivariant cohomology:  $H^*_\Gamma(X) = H^*(X//\Gamma)$.  Moreover $\Phi$ 
is an isomorphism if $\Gamma$ acts freely. Here ${\Cal O}(M)$ denotes
the complexified orientation sheaf on the smooth manifold $M$.

\item{(5)}  (Excision)
Any two-sided ideal $I$ of a complex algebra $A$ gives rise to a 
periodic six-term exact sequence of periodic cyclic cohomology groups
$$
\ldots \to {HC}_{per}^i(A)\to 
{HC}_{per}^i(I){\overset \partial \to \longrightarrow} 
{HC}_{per}^{i+1}(A/I)\to 
{HC}_{per}^{i+1}(A)\to 
\ldots
$$
$i\in \Bbb Z_2$ similar to the topological $K$-theory exact sequence. 
Thus periodic cyclic cohomology defines a generalized cohomology
theory for algebras \cite{CQ,CQ1,CQ2}. This boundary map is multiplicative:
If $B$ is an other algebra and we denote by $\partial_{A\otimes B}$
the boundary map for the exact sequence corresponding to
$I\otimes B \subset A\otimes B$ then we have
$\partial_{A\otimes B}(x\otimes y)=\partial(x) \otimes y$ for any
$x \in HC^*_{per}(I)$ and any $y \in HC^*_{per}(B)$. Here 
$\otimes$ denotes also the external product in cyclic cohomology.

\item{(6)}  There is a functorial morphism 
$$
Ch:  K_i^{\text{alg}}(A)\ \longrightarrow\ HC_i^{per}(A),\quad 
i \in \{ 0,1\}
$$
called the Connes-Karoubi character \cite{C2,K} which is onto for 
$A = {\Bbb C}$  or  $A = {\Bbb C}[{\Bbb Z}]$.  For  $A = C^\infty(X)$, $Ch$ 
coincides with the classical Chern character
up to rescaling \cite{MiS}.

\item{(7)} The boundary map in algebraic  $K$-theory and periodic cyclic 
cohomology are compatible in the following sense
$$
\langle Ind[u],\ \xi\rangle = \langle Ch[u], \partial\xi\rangle
$$
for any $u \in K_1^{alg}(A/I)$  and $\xi \in HC^0_{per}(I)$ \cite{N3}. 
Here we have denoted for simplicity $\langle Ch(Ind[u]),\ \xi\rangle  
=\langle Ind[u],\ \xi\rangle $.

\bigpagebreak

We now go back to our foliation  $(V, {\Cal F})$.  

A complete transversal $N \subset V$  is a submanifold of dimension 
$q=$the codimension of ${\Cal F}$  
which is transverse to the leaves and which intersects each leaf. 
Complete transversals always exist but they are usually not compact and not 
connected.  The choice of a transversal $N$ determines an \'etale groupoid by 
restriction:  ${\Cal G}(N) = {\Cal F}\big\vert_N$.  Explicitly 
${\Cal G}^{(0)}(N) = N$,\ ${\Cal G}^{(1)}(N) = \{ [(x,y,\gamma)] \in 
{\Cal F}$,\ $x,y \in N\}$.  

The equivalence relation is easier to describe in this case.  A path 
$\gamma$  from $x$ to $y$ can be covered by distinguished coordinate patches 
and hence defines a diffeomorphism  $\varphi_\gamma:  N_x \to N_y$  from a 
neighborhood of $x$ in $N$ to a neighborhood of $y$.   Then  
$[(x,y,\gamma)] = [(x_1,y_1,\gamma_1)]$  if and only if  $x = x_1, y = y_1$ 
and  $\varphi_\gamma = \varphi_{\gamma_1}$  in a (possibly smaller) 
neighborhood of $x$.  It is known that the choice of $N$ is not important 
because of Morita equivalence, which gives that  $B{\Cal G}(N)$  is homotopy 
equivalent to  $B{\Cal F}$ \cite{Ha}.

There is a map 
$$
HC^i_{per}(C_c^\infty({\Cal G}(N)) \to HC^i_{per}(\Psi^{-\infty}({\Cal F}))
$$
which is also given by Morita equivalence and an inclusion
$$
j: \Psi^{-\infty}({\Cal F}) \hookrightarrow C_c^\infty({\Cal G}(N)) 
\widehat{\otimes} {\Cal R}
$$
(where ${\Cal R} = \Psi^{-\infty}({\Bbb R}^p),\ p = dim {\Cal F})$  as a full 
corner \cite{BGR}.  

Denote by $q$  the codimension of  ${\Cal F} = dim(N)$.   

\bigpagebreak
\proclaim{Lemma 1}   The morphism
$$
\Phi_0:  H^i(B{\Cal F}, {\Cal O}({\Cal F})) \to HC^{i+q}_{per}(\Psi^{-\infty}({\Cal F}))
$$
given by the composition 
$$
H^i(B{\Cal F}, {\Cal O}({\Cal F})) \simeq H^i(B{\Cal G}(N), {\Cal O}) 
{\buildrel {\Phi}\over{\longrightarrow}}\ 
HC^{i+q}_{per}(C_c^\infty({\Cal G}(N))) 
\to\ HC^{i+q}_{per}(\Psi^{-\infty})({\Cal F}))
$$
does not depend on the choice of $N$.  
\endproclaim

\demo{Proof}  For any complete transversal $N$ the morphism
$$
j: \Psi^{-\infty}({\Cal F}) \hookrightarrow C_c^\infty({\Cal G}(N)) 
\widehat{\otimes} {\Cal R}
$$
depends on a partition of unity in such a way that any two such morphisms 
are conjugated by an inner automorphism.

Let  $N_1$  and $N_2$  be two complete transversals.  Choose a third transversal 
$N'$  not intersecting $N_1$ and $N_2$.  By considering  $N'_1 = N_1 \cup N'$ 
and  $N'_2 = N_2 \cup N'$  one can reduce the problem to the case when  
$N_1 \subset N_2$  and then use the remark in the beginning.\sqep 
\enddemo

\bigpagebreak

Consider the continuous map  ${f}: V \to B {\Cal F}$  which classifies 
${\Cal F}$  \cite{Ha}.  Denote by  $p: S^*{\Cal F} \to V$  the canonical 
projection.  The canonical orientation of the leaves of  $S^*{\Cal F}$  gives 
an identification of  ${\Cal O}(S^*{\Cal F})$  with 
$({f} \circ p)^*{\Cal O}({\Cal F})$.  We use this identification to define a 
morphism
$$
({f} \circ p)^*:  H^*(B{\Cal F}, {\Cal O}({\Cal F})) \to H^*(S^*{\Cal F}, 
{\Cal O}(S^*{\Cal F})).
$$
In the following statement we are going to use the 
notation: 
$$
\Phi_0:  H^i(B{\Cal F}, {\Cal O}({\Cal F})) \to H_{per}^{i+q}(\Psi^{-\infty}
({\Cal F}))
$$ 
for the morphism defined in the previous lemma.  
Also $Ind$ will denote the boundary map in topological 
$K$-theory (as in FOL-TOP), and  
${\Cal T}({\Cal F^*})$ will be the Todd\ class of 
${\Cal F^*}=\text{Hom}_{\Bbb R}({\Cal F},{\Bbb C})$.

\bigpagebreak

\proclaim{Index Formula Problem}  Let ${\Cal F}$  be a foliation
of dimension $n$ and codimension $q$, and let 
$f:V\to B{\Cal F}$, $p:S^*{\Cal F} \to V$ and $\Phi_0$ be as above.
Then for any $u \in K^1(S^*{\Cal F})$ and
$\xi \in H^{2m}(B{\Cal F}, {\Cal O}({\Cal F}))$ we have the following index formula
$$
\langle Ind[u], \Phi_0(\xi)\rangle = {\frac{(-1)^n}{(2\pi i)^m}} 
\langle ch[u] p^*({\Cal T}({\Cal F}^*)f^*(\xi)), [S^*{\Cal F}]\rangle
$$
where $ch$ is the classical Chern character.
\endproclaim

\bigpagebreak

\noindent
The above formula, if correct, would identify the morphism $\ell$  
in Connes' higher index theorem for foliations \cite{C3}.



\bigpagebreak
\subheading{3. The index theorem for foliated bundles}

\bigpagebreak
In case ${\Cal F}$  is actually a fiber bundle the formula in the
above problem becomes the 
formula for the Chern character of the index-bundle of a family.  If $\xi$ 
corresponds to a holonomy invariant measure this 
is Connes' theorem explained in the previous section.  For 
invariant forms in general it is a theorem of Heitsch \cite{He}.  
Previous results for flat bundles were obtained in 
\cite{DHK,J,MN} and for foliations in general but particular 
cocycles in \cite{NT}.  
Except for families of elliptic 
operators, the proofs of all the other circumstances when the index theorem 
for foliations is known are based on heat-kernels.  Our proof marks a 
departure from this approach.  

Consider a smooth fiber bundle $(V, \pi)$,\ $0 \to F \to V\ 
{\buildrel\pi\over\longrightarrow}\ B \to 0$,  with fiber $F$ of dimension 
$q$.  Assume that  $(V,\pi)$  is endowed with a flat connection that we 
interpret as an integrable subbundle  ${\Cal F} \subset TV$, of dimension 
$n = dim(B)$, transverse at any point to the fibers.  The pair $(V,{\Cal F})$  
defines a foliation.  

Fix $b_0 \in B$  and denote by  $\Gamma = \pi_1(B,b_0)$.  Flat fiber bundles 
$(V,\pi,{\Cal F})$  with fiber  $F = \pi^{-1}(b_0)$, $b_0 \in B$, are in a 
one-to-one correspondence with group morphisms 
$\Gamma \to {Diff}(F)$.  Denote 
by  $(\gamma,x) \to \gamma x$ the corresponding action.  Then  $V = 
\widetilde{B} \times_\Gamma F$  where $\Gamma$  acts on the right by deck 
transformation on the universal covering space $\widetilde{B}$  of $B$.

Denote by $\Psi_c^p = \Psi_c^p(\widetilde{B})$  
the space of compactly supported 
order $p$ classical pseudodifferential operators on $\widetilde{B}$.  
Compactly supported 
means that their Schwartz kernel is a compactly supported distribution on 
$\widetilde{B} \times \widetilde{B}$.

Our proof of the index theorem for flat bundles will use 
an alternative description 
of the various algebras associated to the foliation
in terms of certain crossed products.  If the groups $\Gamma$ acts on an 
algebra  $A_0$  then  the {\it algebraic crossed product} 
$A_0 \rtimes \Gamma$  is by definition the linear span of formal 
products  $\{ a\gamma,\ a \in A_0,\ \gamma \in \Gamma\}$  with the product 
rule  $(a\gamma)(a'\gamma') = a\gamma(a') \gamma\gamma'$.  

The group $\Gamma$  acts on the spaces  
$\Psi_c^p = \Psi_c^p(\widetilde{B})$ 
so we obtain exact sequences

$$0 \to \Psi_c^{-\infty} \rtimes \Gamma \to 
\Psi_c^0 \rtimes \Gamma \to (\Psi_c^0/\Psi_c^{-\infty})\rtimes \Gamma \to 0$$

$$0 \to \Psi_c^{-1} \rtimes \Gamma \to \Psi_c^0 \rtimes \Gamma\ 
{\buildrel {\sigma_0}\over\longrightarrow}\ C_c^\infty(S^*\widetilde{B})\ 
\rtimes \Gamma \to 0$$

\noindent
where $\sigma_0$  is the principal symbol map and  $S^*\widetilde{B} \subset 
T^*\widetilde{B}$  is the cosphere bundle of $\widetilde{B}$.  

By a standard procedure we enlarge the algebra  $\Psi_c^0$  to include all 
$(n+1)$-summable Schatten-von Neumann operators.  Explicitly, denote by 
$$
C_{n+1} = C_{n+1}(L^2(\widetilde{B})) = \{ T, tr(T^*T)^{(n+1)/2} < \infty\}
$$
where  $T$  denotes a bounded operator on  $L^2(\widetilde{B})$
and $tr$ is the usual trace.

It is a simple known fact that  $\Psi_c^{-1} \subset C_{n+1}$  and that 
$\Psi^0_c \cap C_{n+1} \subset \Psi_c^{-1}$.  This tells us that if we define 
$E = E(\widetilde{B})$  we have an exact sequence
$$
0 \to C_{n+1} \to E\ {\buildrel {\sigma_0}\over\longrightarrow}\
C_0^\infty(S^*\widetilde{B}) \to 0
$$
considered also in \cite{N2,N3}.

Fix $2m \ge n+1$  and denote by  $Tr_m \in HC^{2m}(C_{n+1})$  the continuous 
cocycle 
$$
Tr_m(a_0,\cdots,a_{2m}) = {\frac{(-1)^mm!}{(2m)!}} tr(a_0 \cdots a_m).
$$
The normalization factor is chosen such that  $Tr_m\big\vert_{C_1} = 
S^mtr$, where  $S$ is the Connes periodicity operator.  By abuse of notation 
we denote by  $Tr \in HC^0_{per}(C_{m+1})$  the class of  $Tr_m$  for any 
$2m \ge n+1$.         

Denote by  $H^*(X)$  the ${\Bbb Z}/2{\Bbb Z}$-periodic 
complex cohomology groups of the a manifold $X$.

\bigpagebreak

\proclaim{Lemma 2}  We have that  $\Phi$  induces an isomorphism 
$$
H^{*-1}(S^*B) \simeq HC^*_{per}(C_c^\infty(S^*\widetilde{B}) \rtimes \Gamma) 
,\ i\in {\Bbb Z}/2{\Bbb Z}.
$$
\endproclaim

\demo{Proof}  Since  $\Gamma$  acts freely on the oriented
odd dimensional manifold 
$S^*\widetilde{B}$  and $(S^*\widetilde{B})/\Gamma = S^*B$  then we have that 
the map  $\Phi$  is an isomorphism, see \cite{BN}.\sqep \enddemo

\bigpagebreak

Recall \cite{N3} that there is a 
$HC^*_{per}({\Bbb C}[\Gamma])$  module structure 
on  $HC^*_{per}(A \rtimes \Gamma)$  induced by the 
${\Bbb C}[\Gamma]$-coalgebra 
structure of   $A \rtimes \Gamma$:
$$
A \rtimes \Gamma \ni a\gamma {\buildrel 
{\delta}\over\longrightarrow}\ a\gamma 
\otimes \gamma \in (A \rtimes \Gamma) 
\otimes {\Bbb C}[\Gamma].
$$
thus 
$xy=\delta^*(x \otimes y)$ if $y\in HC^*_{per}({\Bbb C}[\Gamma])$
and $x\in HC^*_{per}(A \rtimes \Gamma)$.

Denote by  $g: S^*B \to B\Gamma$  the classifying map of the covering 
$\Gamma \to S^*\widetilde{B} \to S^*B$.  Also recall \cite{C2}
that  $H^*(B\Gamma)$  is a 
direct summand of  $HC^*_{per}({\Bbb C}[\Gamma])$.  Denote by 
$$
r_0:  HC^*_{per}({\Bbb C}[\Gamma]) \longrightarrow H^*(B\Gamma)
$$
the natural projection.  It is a ring morphism satisfying
$r_0\circ \Phi = id$ where $\Phi:H^*(B\Gamma)\to 
HC^*_{per}({\Bbb C}[\Gamma])$ is the geometric map.

\bigpagebreak

\proclaim{Lemma 3}  The  $HC^*_{per}({\Bbb C}[\Gamma])$-module structure of 
$HC^*_{per}(C_c^\infty(S^*\widetilde{B}) \rtimes \Gamma)$  is given, using the 
isomorphism of the previous lemma, by  $g^* \circ r_0$:
$
\Phi(\zeta) \xi = \Phi(\zeta g^*(\cup r_0(\xi)))
$
for any  $\zeta \in H^{*-1}(S^*B)$  and  
$\xi \in HC^*_{per}({\Bbb C}[\Gamma])$. In particular 
$\Phi(\zeta)\Phi(\eta)=\Phi(\zeta g^*(\eta))$.
\endproclaim

\demo{Proof}  We know that the action of 
$HC^*_{per}({\Bbb C}[\Gamma])$  factors
through  $r_0$  because $\Gamma$  acts without fixed points on  
$S^*\widetilde{B}$. This shows that we can assume $\xi =\Phi(\eta)$ for 
$\eta=r_0(\xi)$, $\eta\in H^*(B\Gamma)$.  
The module structure is obtained using the multiplicativity of $\Phi$ and 
the fact that the composition 
$$
S^*B = S^*\widetilde{B}//\Gamma\ \longrightarrow\ 
S^*\widetilde{B}//(\Gamma \times \Gamma) = S^*B \times B\Gamma
$$
corresponding to $\delta$ is  $id \times g$, by definition.
We then have 
$$
\Phi(\zeta)\Phi(\eta)=
\delta^*(\Phi(\zeta)\otimes \Phi(\eta))=
\Phi(\zeta g^*(\eta))\sqep 
$$
\enddemo

\bigpagebreak
 
\proclaim{Proposition}  Suppose the graph of $(V, {\Cal F})$
is separated where $V,F,B$ and $E=E(\widetilde{B})$ are as above. Then 
there exists a commutative diagram
$$
\matrix \format \c &\c &\c &\c &\c &\c &\c 
&\c &\c \\
0 & \longrightarrow & \Psi^{-1}({\Cal F}) & \longrightarrow & \Psi^0({\Cal F}) 
& \longrightarrow & C_c^\infty(S^*{\Cal F}) & \longrightarrow & 0 \\
\vspace{2\jot}
&& \alpha_0 \bigg\downarrow && \alpha \bigg\downarrow 
&& \bigg\downarrow \beta \\
\vspace{2\jot}
0 & \longrightarrow & C_{n+1} \widehat{\otimes} C_c^\infty(F)\rtimes \Gamma 
& \longrightarrow & (E \widehat{\otimes} C_0^\infty (F)) \rtimes \Gamma 
& \longrightarrow & C_c^\infty(S^*\widetilde{B} \times F) \rtimes \Gamma 
& \longrightarrow & 0  
\endmatrix
$$
where $\alpha$ is an isomorphism onto 
$e((E \widehat{\otimes} C_0^\infty (F)) \rtimes \Gamma) e$ for some
idempotent $e$ and 
$\beta$ induces an isomorphism in cyclic cohomology.
\endproclaim

\demo{Proof}  This is the promised 
equivalent definition of the various algebras 
associated to  
$({V}, {\Cal F})$  in the particular case of a foliated 
bundles.  The idempotent $e$ is defined using a standard argument
based on partitions of unity as follows.

Choose a partition
of unity on $B$ subordinated to a finite trivializing cover $(U_i)_{i=1,N}$
of $\widetilde{B}$. We can find smooth real functions $\varphi_i$,
$\sum_i \varphi_i^2=1$, 
with support in $U_i$.
Choose {\it disjoint}
open sets $V_i \subset \widetilde{B}$ on which the projection
$\widetilde{B}\to B$ is a
diffeomorphism. Denote by $\gamma_{i, j}$ the corresponding 
2-cocycle with values in $\Gamma$ and by 
$\widetilde{\varphi}_i$ the lifts of $\varphi_i$
to functions supported on $V_i$, so that 
$\widetilde{\varphi}_i\widetilde{\varphi}_j=0$ if $i\not =j$. Then 
$e=\sum_{i,j}\widetilde{\varphi}_i \gamma_{i,j} \widetilde{\varphi}_j$
is the idempotent $e$ we were looking for. Note that if
$B$ is not compact $\gamma_{i,j}$ are going to be locally constant functions
on $U_i\cap U_j$, not constant in general, so $e$ will have  entries in
the multiplier algebra of $E(\widetilde{B})$
and not in $E(\widetilde{B})$ itself.
\sqep \enddemo

\bigpagebreak   

Consider now the exact sequence
$$
0 \to C_{n+1} \rtimes \Gamma \to E \rtimes \Gamma \to C_c^\infty(S^*\widetilde{B}) 
\rtimes \Gamma \to 0.
$$

We want to identify the Cuntz-Quillen boundary map
$$
\partial_{E\rtimes \Gamma}: HC^*_{per}(C_{n+1} \rtimes \Gamma) \to 
HC^{*+1}_{per}(C_c^\infty(S^*\widetilde{B}) \rtimes \Gamma) \simeq H^*(S^*B)
$$
of this exact sequence.

Observe that since $\Gamma$ acts by inner automorphisms on $C_{n+1}$  we have 
that  $C_{n+1} \rtimes \Gamma \simeq C_{n+1} \otimes {\Bbb C}[\Gamma]$.

Define the {\it Index characteristic class}
${\Cal I}(M) \in H^{even}(S^*M)$ of a smooth manifold $M$
to be the rescaled Todd class of
$T^*M$, ${\Cal I}(M)_{2k}
=(2\pi \imath)^{n-k} p^*{\Cal T}(T^*M\otimes \Bbb C)_{2k}$ 
where ${\Cal T}$ is  the usual Todd class and $n=dim(M)$.

\bigpagebreak

\proclaim{Lemma 4}  Let $\xi \in H^j(B\Gamma)$, $n=dim(B)$,  then
$$
\partial_{E\rtimes \Gamma}(Tr \otimes \Phi(\xi)) = (-1)^n\Phi({\Cal I}(B)
g^*(\xi)) \in HC^{j+1}(C_c^\infty(S^*B))
$$
where  ${\Cal I}$  is the Index class and $\Phi:H^{*}(S^*B)
\to HC^{*+1}(C_c^\infty(S^*B))$ is as in the previous section.
\endproclaim

\demo{Proof} 
Consider the following commutative diagram
$$
\matrix \format \c &\quad \c &\quad \c &\quad \c &\quad \c &\quad \c &\quad \c 
&\quad \c &\quad \c \\
0 & \longrightarrow & C_{n+1}  & \longrightarrow & E(B) 
& \longrightarrow & C_c^\infty(S^*B) & \longrightarrow & 0 \\
\vspace{2\jot}
&& \alpha_0' \bigg\uparrow && \alpha' \bigg\uparrow && \bigg\uparrow = \\
\vspace{2\jot}
0 & \longrightarrow & e(C_{n+1}\rtimes \Gamma)e  & \longrightarrow & {\frak A} 
& \longrightarrow & C_c^\infty(S^*B) & \longrightarrow & 0 \\
\vspace{2\jot}
&& \alpha_0 \bigg\downarrow && \alpha \bigg\downarrow && \bigg\downarrow \beta \\
\vspace{2\jot}
0 & \longrightarrow & C_{n+1} \rtimes \Gamma & \longrightarrow & 
 E(\widetilde{B}) \rtimes \Gamma & \longrightarrow & 
 C_c^\infty(S^*\widetilde{B}) \rtimes \Gamma & \longrightarrow & 0  
\endmatrix
$$
where ${\frak A}=e(E(\widetilde{B}) \rtimes \Gamma)e$,
for an idempotent $e$ implementing the Morita equivalence, 
and $\alpha$ is the inclusion. The idempotent $e$ is defined
as in the previous proposition. The morphism $\alpha'$
is defined using the natural
representation of $E(\widetilde{B}) \rtimes \Gamma$
on $L^2(\widetilde{B})$ given by the fact that the action
of $\Gamma$ is implemented by inner automorphisms. This shows that
the restriction $\alpha_0'$ is given by the composition
$$
e(C_{n+1}\rtimes \Gamma)e \to C_{n+1}\rtimes \Gamma \simeq 
C_{n+1}\otimes \Bbb C[\Gamma] {\buildrel {1 \otimes \chi}\over{\longrightarrow}}
C_{n+1}
$$
where $\chi : \Bbb C[\Gamma] \to \Bbb C$  is the augmentation morphism
$\gamma \to 1$.

The above  commutative diagram has the property that
$$
\beta^*:  HC^*_{per}(C_c^\infty(S^*\widetilde{B}) \rtimes \Gamma) 
\longrightarrow
HC^*_{per}(C_c^\infty(S^*B))
$$
is the isomorphism of Lemma 2 and that $\alpha_0^*$ is an isomorphism
as well. Also, if we regard the the augmentation morphism 
$\chi : \Bbb C[\Gamma] \to \Bbb C$ 
as a trace (an element of cyclic cohomology) then
$\alpha_0^*(Tr \otimes \chi) = \alpha_0^{'*}(Tr)$, and  
$\chi$ is the identity of $HC^*_{per}(\Bbb C[\Gamma])$. This
gives that $(Tr \otimes \chi) \Phi(\xi)=Tr \otimes \Phi(\xi)$.
The boundary $\partial_{AS}$  of the top exact sequence is determined by the 
Atiyah-Singer formula which gives, using the
results of \cite{N3}, $\partial_{AS}(Tr)=(-1)^n\Phi({\Cal I}(B))$.
Denote by $\partial_2$
the boundary map of the middle exact sequence. Using the
$HC^*_{per}(\Bbb C[\Gamma])$-linearity of $\partial_2$ and 
$\partial_{E \rtimes \Gamma}$, see \cite{Ni3}, we obtain
$$
\beta^{*}\circ \partial_{E \rtimes \Gamma}(Tr \otimes \Phi(\xi)) =
\partial_{2} 
\circ \alpha_0^*(Tr \otimes \chi)\Phi(\xi)= \partial_{2}\circ\alpha_0^{'*}(Tr)\Phi(\xi)
=\partial_{AS}(Tr)\Phi(\xi)
$$
Since  $\partial_{AS}(Tr)\Phi(\xi)=(-1)^n\Phi({\Cal I}(B)) \Phi(\xi)=
(-1)^n\Phi({\Cal I}(B) g^*(\xi))$, using the formula from the previous lemma,
the result follows.\sqep \enddemo

\bigpagebreak

It is interesting to observe that the mere existence of the top
commutative diagram in the previous proof implies a theorem of
Atiyah and Singer \cite{A}. The lemma is equivalent to the higher
index theorem for coverings of Connes and Moscovici \cite{CM1,N3}.

Suppose now that the foliation  $({V}, {\Cal F})$ defined 
at the beginning of this section by a flat connection of the
bundle $\pi:V\to B$ has a separated graph. 
This is equivalent to the following condition:  the only $\gamma \in \Gamma$  
for which $F^\gamma$  has a nonempty interior are those $\gamma$ that act 
trivially on $F$.  Here $F$ denotes  the fiber of $V \to B$ as before.

\bigpagebreak

Consider now the exact sequence 
$$
0 \to C_{n+1} \widehat{\otimes} C_c^\infty(F) \rtimes \Gamma \to E({\Cal F})\ 
{\buildrel {\sigma_0}\over\longrightarrow}\ C_c^\infty(S^*{\Cal F}) \to 0
$$
induced by the morphism $\beta$ in the previous proposition.  Denote by
$$
\partial_{\Cal F}: HC^*_{per}(C_{n+1} \widehat{\otimes} C_c^\infty(F) \rtimes 
\Gamma) \longrightarrow HC_{per}^{*+1}(C_c^\infty(S^*{\Cal F})) = H^{*+q}(S^*{\Cal F}, 
{\Cal O}(S^*{\Cal F}))
$$
where $q$  is the codimension of ${\Cal F}$, and ${\Cal O}(S^*{\Cal F})$  is the 
orientation sheaf of  $S^*{\Cal F}$.  

In the particular case of $(V,{\Cal F})$ that we are studying, the classifying
space $B{\Cal F}$ 
of the graph of the foliation coincides, up to homotopy, with
the homotopy quotient $F//\Gamma$. The map $f:V\to B{\Cal F}$ can be
described as the second component in
the composition  $(id,f):V = (\widetilde{B} \times F)/\Gamma \to 
(\widetilde{B} \times F)// (\Gamma \times \Gamma) = B \times (F // \Gamma)$. 
This map satisfies 
${\Cal O}(S^*{\Cal F}) \simeq p^*f^* {\Cal O}(F)$  
where the isomorphism depends on a 
choice of the orientation of the leaves of  $S^*{\Cal F}$, and  
$p:S^*{\Cal F}\to V$  is the projection. We orient the leaves of 
$S^*{\Cal F}$ as the boundary of an almost complex manifold.

Also recall from the previous section that there 
exist maps  
$$
\Phi:  H^*(F//\Gamma, {\Cal O}(F)) \to 
HC^{*+q}_{per}(C_c^\infty(F) \rtimes \Gamma)
$$ 
and
$$
\Phi_0:  H^*(B{\Cal F},{\Cal O}({\Cal F})) \to 
HC^{*+q}_{per}(\Psi^{-\infty}({\Cal F}))
$$ 
where $q=dim (F)$.
In the case we discuss now, that of a foliated bundle, these
two maps are related by $Tr \otimes \Phi(\xi) =\Phi_0(\xi)$.

Denote by $\pi_0:S^*{\Cal F}\to S^*B$ the projection covering
$\pi:V\to B$ and let ${\Cal I}({\Cal F})=\pi_0^*({\Cal I}(B))\in 
H^*(S^*{\Cal F})$ be the Index class of ${\Cal F}$.

\bigpagebreak
\proclaim{Theorem}  Let  $\xi \in H^*(F//\Gamma, {\Cal O}(F))$  and 
$u \in K_1(C^\infty(S^*{\Cal F}))$  be an invertible symbol, also let 
$f:V\to B{\Cal F}=F//\Gamma$ and $p:S^*{\Cal F} \to V$ be 
as above, $f_0=f\circ p$.  
Then 
$$
\langle \Phi_0(\xi),\ Ind[u]\rangle = (-1)^n
\langle Ch[u]{\Cal I}({\Cal F}) 
f_0^*(\xi)), [S^*{\Cal F}]\rangle
$$
where  $Ind: K_1^{alg}(C^\infty(S^*{\Cal F})) 
\to K_0^{alg}(C_{n+1} \otimes C^\infty(F) 
\rtimes \Gamma)$  is the boundary map in algebraic $K$-theory,
$n=dim (B)$ and the cup product is defined using 
${\Cal O}(S^*{\Cal F}) \simeq f^*_0 {\Cal O}(F)$.  
\endproclaim

\demo{Proof}  Denote $\Gamma_2=\Gamma \times \Gamma$.
Consider the following commutative diagram:
$$
\matrix \format \c &\c &\c &\c &\c &\c &\c 
&\c &\c \\
0 & \longrightarrow & C_{n+1} \otimes C_c^\infty(F) \rtimes \Gamma  
& \longrightarrow & (E \otimes C_c^\infty(F)) \rtimes \Gamma
& \longrightarrow & C_c^\infty(S^*\widetilde{B} \times F)\rtimes \Gamma 
& \longrightarrow & 0 \\
\vspace{2\jot}
&& \alpha_0 \bigg\downarrow && \alpha \bigg\downarrow && \bigg\downarrow \beta \\
\vspace{2\jot}
0 & \longrightarrow & (C_{n+1} \otimes C_c^\infty(F)) \rtimes \Gamma_2
& \longrightarrow & (E \otimes C_c^\infty (F)) \rtimes \Gamma_2 	 
& \longrightarrow & C_c^\infty(S^*\widetilde{B} \times F) \rtimes 
\Gamma_2 & \longrightarrow & 0  
\endmatrix
$$	       

where the vertical arrows are defined 
by the diagonal morphism  $\Gamma \to \Gamma_2=\Gamma \times \Gamma$
and
the bottom line is obtained from the exact sequence of  
$E(\widetilde{B}) \rtimes \Gamma$ by 
tensoring with  $C_c^\infty(F) \rtimes \Gamma$.  

The map $(S^*\widetilde{B}\times F)//\Gamma \to 
(S^*\widetilde{B}\times F)//\Gamma_2$ is $(\pi_0,f_0):
S^*{\Cal F}\to S^*B\times (F//\Gamma)$ which shows that
the morphism $\beta^*$  satisfies  the relation
$$
\beta^*(\Phi(\xi_0) \otimes \Phi(\xi_1))
=\Phi(\pi_0^*(\xi_0) f_0^*(\xi_1)) \in 
HC_{per}^*(C_c^\infty(S^*\widetilde{B} \times F)\rtimes \Gamma) 
$$
where 
$\xi_0 \in H^{*+1}(S^*B)$  and $\xi_1\in H^{*+q}(F//\Gamma, 
{\Cal O}(F))$. 
This follows from the functoriality and 
multiplicativity of the geometric morphism $\Phi$.

We have, using the multiplicativity of
the boundary morphism and Lemma 4,

$$
\aligned
\partial_{\Cal F}(\Phi_0(\xi)) &=\partial_{\Cal F}(Tr \otimes \Phi(\xi)) \\
  &= \partial_{\Cal F}(\alpha_0^*(Tr 
\otimes 1 \otimes \Phi(\xi))) \\
  &= \beta^*(\partial_{E\rtimes \Gamma}(Tr \otimes 1) \otimes \Phi(\xi)) \\
  &= (-1)^n\beta^*(\Phi({\Cal I}(B) g^*(1))\otimes \Phi(\xi))\\ 
  &= (-1)^n\Phi(\pi_0^*({\Cal I}(B)) f_0^*(\xi))\\
  &= (-1)^n\Phi({\Cal I}({\Cal F}) f_0^*(\xi)).
\endaligned
$$

(here $g:S^*B\to B\Gamma$ classifies the universal covering of $S^*B$
as in Lemma 4).

We obtain, using the compatibility
between the index map in $K$-Theory and the boundary map in 
periodic cyclic cohomology,

$$
\aligned
\langle \Phi_0(\xi), Ind[u]\rangle &= 
\langle \partial_{\Cal F}(\Phi_0(\xi)), Ch[u]\rangle\\
&= (-1)^n\langle \Phi({\Cal I}({\Cal F}) f_0^*(\xi)), Ch[u]\rangle \\
&= (-1)^n\langle Ch[u] {\Cal I}({\Cal F}) f_0^*(\xi), 
[S^*{\Cal F}]\rangle .
\endaligned
$$

This proves the theorem.\sqep 
\enddemo



\widestnumber\key{DHK}
\Refs
\ref\key A 
          \by {M. F. Atiyah}
          \paper {The signature of fiber bundles}
          \jour {Global Analysis, papers in honour of K. Kodaira},
          {Univeristy of Tokio Press and Princeton University
         Press}
          \yr {1969}
          \pages {73-84}
\endref
\ref\key AS  
          \by {M. F. Atiyah and I. M. Singer}
          \paper {The index of elliptic operators {III}}
          \jour {Ann. of Math.}
          \yr {1968} 
          \vol {87}
          \pages {546-604}
\endref 
\ref\key BGR %
          \by {L. Brown, P. Green, M. Rieffel}
          \paper {Stable isomorphism and strong Morita 
equivalence of $C^*$-algebras}
          \vol {71}
          \jour {Pacific J. Math.}
          \pages {349-363}
\endref
\ref\key BN  
          \by {J.-{L}. Brylinski and V. Nistor}
          \paper {Cyclic cohomology of etale groupoids}
          \vol {8}
          \jour {K-{T}heory}
          \pages {341-365}
          \yr {1994}
\endref  
\ref\key C1 
          \by {A. Connes}
          \paper {Sur la th\'{e}orie noncommutative de l'int\'{e}gration}
          \jour {Lect. Notes Math.}
          \yr {1982}
          \vol {725}
          \pages {19-143}
\endref 
\ref\key C2 
          \by {A. Connes}
          \paper {Non-commutative differential geometry}
          \jour {Publ. Math. IHES}
          \yr {1985}
          \vol {62}
          \pages {41-144}
\endref
\ref\key C3 
          \by {A. Connes}
          \paper {Non {C}ommutative {G}eometry}
          \jour {Academic Press}
          \yr{1995}
\endref 
\ref\key C4 
          \by {A. Connes}
          \paper {Noncommutative geometry and reality}
          \jour {preprint IHES}
          \yr{1995}
\endref
\ref\key CM1 
          \by {A. Connes and H. Moscovici}
          \paper {Cyclic cohomology, the Novikov conjecture and hyperbolic groups}
          \jour {Topology}
          \vol {29}
          \yr{1990}
          \pages{345-388}
\endref
\ref\key CM2 
          \by {A. Connes and H. Moscovici}
          \paper {The local index formula in noncommutative geometry}
          \jour {Geom. and Funct. Anal.}
          \yr{1995}
          \pages{170-237}
\endref
\ref\key CQ
          \by {J. Cuntz and D. Quillen}
          \paper {On excision in periodic cyclic cohomology}
          \jour {preprint}
          \yr {1994}
\endref 
\ref\key CQ1
          \by {J. Cuntz and D. Quillen}
          \paper {On excision in periodic cyclic cohomology}
          \jour {C. R. Acad. Sci. Paris}
          \yr {1993}
          \vol {317}
          \pages {917-922}
\endref 
\ref\key CQ2
          \by {J. Cuntz and D. Quillen}
          \paper {On excision in periodic cyclic cohomology {II}, the
		  general case}
          \jour {C. R. Acad. Sci. Paris}
          \yr {1994}
          \vol {318}
          \pages {11-12}
\endref
\ref\key DHK
          \by {R. Douglas, S. Hurder, J. Kaminker}
          \paper {Cyclic
cocycles, renormalization
and eta-invariants}
          \jour {Invent. Math.}
          \yr {1991}
          \vol {103}
          \pages {101-179}
\endref
\ref\key Ha 
          \by {A. Haefliger}
          \paper {Groupo\"{i}des d'holonomie et espaces classifiants}
          \jour {Ast\'{e}risque}
          \yr {1984}
          \vol {116}
          \pages {70-97}
\endref 
\ref\key He %
          \by {J. Heitsch}
          \paper {Bismut superconnections and the
Chern character for Dirac operators on foliated manifolds}
          \jour {Preprint}
\endref
\ref\key J
          \by {X. Jiang}
          \paper {An index theorem for foliated bundles}
          \jour {preprint}
\endref
\ref\key K %
          \by {M. Karoubi}
          \paper
{Homologie cyclique  et K-theorie}
          \vol 149
          \jour {Asterisque}
          \pages {1-147}
          \yr 1987
\endref  
\ref\key LQ
          \by {{J.-L.} Loday and D. Quillen}
          \paper {Cyclic homology and the {L}ie algebra homology
                     of matrices}
          \jour {Coment. Math. Helveticii}
          \yr {1984}
          \vol {59}
          \pages {565-591}
\endref 
\ref\key MiS 
          \by {J. Milnor and J. Stasheff}
          \paper {Characteristic Classes}
          \yr {1974}
          \jour {Ann. Math. Studies}
          {Princeton University Press}
          \vol {76}
          \endref
\ref\key MoS%
          \by {C. C. Moore, C. Schochet}
          \paper {Global Analysis on Foliated Spaces}
          \jour {Math. Sci. Res. Inst. Series}, {Springer Verlag}
          \vol{ }
\endref
\ref\key MN%
          \by {H. Moriyoshi, T. Natsume}
          \paper {The Godbillion-Vey cyclic cocycle and longitudinal
          Dirac operators}
          \jour {to appear in Pacific J. Math.}
\endref
\ref\key  NT%
          \by {R. Nest, B. Tsygan}
          \paper
{Algebraic index theorem
for families}
          \jour {Adv.
Math. (to appear)}
\endref  
\ref\key N1 
          \by {V.  Nistor}
          \paper {Group cohomology and the cyclic cohomology of crossed
                     products}
          \jour {Invent. Math.}
          \yr {1990}
          \vol {99}
          \pages {411-424}
\endref 
\ref\key N2 
          \by {V. Nistor}
          \paper {A bivariant {Chern-Connes} character}
          \yr {1993}
          \vol {138}
          \jour {Ann. of Math.}
          \pages {555-590}
\endref 
\ref\key N3 
          \by {V. Nistor}
          \paper {Properties of the Cuntz-Quillen boundary map}
          \jour {preprint}
\endref
\ref\key W %
          \by {E. Wilkenkemper}
          \paper {The graph of a foliation}
          \vol {1}
          \jour {Ann. Global Anal. Geom.}
\pages 53-73
\endref
\endRefs
\enddocument